\documentclass{amsart}

\usepackage{fancyhdr}
\usepackage{lastpage}
\usepackage{stmaryrd,yhmath}

\newcommand{\bdis}{\begin{displaymath}}
\newcommand{\edis}{\end{displaymath}}
\newcommand{\be}{\begin{equation}}
\newcommand{\ee}{\end{equation}}
\newcommand{\mbb}{\mathbb}
\newcommand{\mcal}{\mathcal}

\newcommand{\vth}{\vartheta}

\newcommand{\zf}{\zeta\left(\frac{1}{2}+it\right)}

\theoremstyle{definition}

\theoremstyle{remark}
\newtheorem{remark}[]{Remark}

\newtheorem*{mydef1}{{\bf Theorem}}

\newtheorem*{mydef41}{{\bf Corollary 1}}

\newtheorem*{mydef42}{{\bf Corollary 2}}

\newtheorem*{mydef43}{{\bf Corollary 3}}

\newtheorem*{mydef44}{{\bf Corollary 4}}

\numberwithin{equation}{section}

\begin{document}

\title{Jacob's ladders, new properties of the function $\arg\zf$ and corresponding metamorphoses}

\author{Jan Moser}

\address{Department of Mathematical Analysis and Numerical Mathematics, Comenius University, Mlynska Dolina M105, 842 48 Bratislava, SLOVAKIA}

\email{jan.mozer@fmph.uniba.sk}


\begin{abstract}
The notion of the Jacob's ladders, reversely iterated integrals and the $\zeta$-factorization is used in this paper in order to obtain new results
in study of the function $\arg\zf$. Namely, we obtain new formulae for non-local and non-linear interaction of the functions
$|\zf|$ and $\arg\zf$, and also a set of metamorphoses of the oscillating Q-system.
\end{abstract}
\maketitle

\section{Introduction}

\subsection{}

Let us denote by $N(T)$ the number of zeroes $\beta+i\gamma$ of the $\zeta(s)$-function such that
\bdis
\beta\in (0,1),\ \gamma\in (0,T).
\edis
We suppose that $T$ is not equal to any $\gamma$. Otherwise, we put
\bdis
N(T)=\frac 12\lim_{\epsilon\to 0^+} [N(T+\epsilon)+N(T-\epsilon)].
\edis
It us well-know that
\bdis
N(T)=\frac{1}{2\pi}T\ln\frac{T}{2\pi e}+\frac 78+S(T)+\mcal{O}\left(\frac 1T\right),
\edis
where
\be \label{1.1}
S(T)=\frac 1\pi \arg\zf,
\ee
and the value of $\arg$ is obtained by continuous variation along the straight lines joining the points
\bdis
2,\ 2+iT,\ \frac 12+iT ,
\edis
starting with the value zero. Next, we have the function
\be \label{1.2}
S_1(T)=\int_0^T S(t){\rm d}t.
\ee

\subsection{}

Further, let us remind the following facts
\bdis
\zf=\left|\zf\right|e^{i\arg\zf},
\edis
i.e. the functions
\be \label{1.3}
\left|\zf\right|,\ \arg\zf
\ee
are parts of the Riemann function
\bdis
\zf.
\edis
The study of these functions have proceeded by isolated ways. Namely:
\begin{itemize}
 \item[(a)] the first one studied by Hardy-Littlewood
 \bdis
 \int_T^{T+U} \left|\zf\right|^2{\rm d}t\sim U\ln T, \dots
 \edis
 \item[(b)] the second one by Backlund, E. Landau, H. Bohr, Littlewood, Titchmarsh - to the fundamental
 Selberg's results (see \cite{3}, \cite{4}).
\end{itemize}

Let us mention two of the Selberg's results:
\be \label{1.4}
\begin{split}
 & \int_T^{T+H}\{ S_1(t)\}^{2l}{\rm d}t=c_lH+\mcal{O}\left(\frac{H}{\ln T}\right), \\
 & T^a\leq H\leq T;\ \frac 12<a\leq 1,\ l\in\mbb{N},
\end{split}
\ee
where $l$ is arbitrary and fixed, (see \cite{4}, p. 130), and
\be \label{1.5}
S_1(t)=\Omega_{\pm}\left\{(\ln t)^{1/3}(\ln\ln t)^{-10/3}\right\},
\ee
(see \cite{4}, p.150).

\begin{remark}
For our purpose it is sufficient to use the formula (\ref{1.4}) in the minimal case
\bdis
H=T^{1/2+\epsilon},\ \epsilon>0,
\edis
where $\epsilon$ is sufficiently small (non-principal improvements of the exponent $1/2$ are not
relevant for our purpose).
\end{remark}

\subsection{}

To this date, there is no result in the theory of the Riemann zeta-function about the interaction of the
functions (\ref{1.3}), or the functions
\be \label{1.6}
\left|\zf\right|, S_1(t).
\ee
That is, there is nothing known like
\bdis
F\left(\left|\zf\right|,|S_1(\tau)|\right)=0
\edis
for a set of values $t,\tau$.

On the other hand, we have developed (see \cite{1}) the method of $\zeta$-factorization that gives, for example,
the following formula (see \cite{1}, (1.7))
\bdis
\frac{1}{\sqrt{\left|\zeta\left(\frac 12+i\alpha_0\right)\right|}}\sim\frac{1}{\sqrt{\Lambda}}\prod_{r=1}^k
\left|\zeta\left(\frac 12+i\alpha_r\right)\right|
\edis
together with the infinite set of corresponding metamorphoses of the main multiform.

In this paper we use this method to obtain a result of the new type
\be \label{1.7}
|S_1(\alpha_0)|\sim \Phi
\left\{
\prod_{r=1}^k
\left|
\frac
{\zeta\left(\frac 12+i\alpha_r\right)}
{\zeta\left(\frac 12+i\beta_r\right)}
\right|
\right\}
\ee
together with the infinite set of metamorphoses of the corresponding Q-system from \cite{2}.

\begin{remark}
A kind of nonlinear and nonlocal interaction of the functions (\ref{1.6}) is expressed by the formula (\ref{1.7}).
\end{remark}

\section{Theorem}

\subsection{}

We begin with the Selberg's formula
\be \label{2.1}
\begin{split}
 & \int_T^{T+H}\{ S_1(t)\}^{2l}{\rm d}t\sim c_lH,\ T\to\infty , \\
 & H=T^{1/2+\epsilon},\ l\in\mbb{N},\ \epsilon>0,
\end{split}
\ee
(comp. (\ref{1.4}) and Remark 1), where $l$ is arbitrary and fixed, $\epsilon$ is sufficiently small. Now, if we use
our method of transformation (see \cite{2}, (4.1)--(4.19)) in the case of the formula (\ref{2.1}) then we obtain
(see (\ref{1.1}), (\ref{1.2})) the following

\begin{mydef1}
Let
\be \label{2.2}
[T,T+H]\longrightarrow [\overset{1}{T},\overset{1}{\wideparen{T+H}}],\dots , [\overset{k}{T},\overset{k}{\wideparen{T+H}}],
\ee
where
\bdis
[\overset{r}{T},\overset{r}{\wideparen{T+H}}],\ r=1,\dots,k,\ k\leq k_0,\ k_0\in\mbb{N}
\edis
be the reversely iterated segment corresponding to the first segment in (\ref{2.2}) and $k_0$ be an arbitrary and fixed
number. Then there is a sufficiently big
\bdis
T_0=T_0(l,\epsilon)>0
\edis
such that for every $T>T_0$ and every admissible $l,\epsilon,k$ there are the functions
\be \label{2.3}
\begin{split}
 & \alpha_r=\alpha_r(T,l;\epsilon,k),\ r=0,1,\dots,k, \\
 & \beta_r=\beta_r(T;\epsilon,k),\ r=1,\dots,k, \\
 & \alpha_r,\beta_r\not=\gamma:\ \zeta\left(\frac 12+i\gamma\right)=0
\end{split}
\ee
such that
\be \label{2.4}
\begin{split}
 & \left|\int_0^{\alpha_0(T)} \arg\zf {\rm d}t\right|\sim \\
 & \sim \pi (c_l)^{\frac{1}{2l}}\prod_{r=1}^k
 \left|
\frac
{\zeta\left(\frac 12+i\alpha_r(T,l)\right)}
{\zeta\left(\frac 12+i\beta_r(T)\right)}
\right|^{-\frac 1l},\ T\to\infty .
\end{split}
\ee
Moreover, the sequences
\bdis
\{ \alpha_r\}_{r=0}^k,\ \{\beta_r\}_{r=1}^k
\edis
have the following properties
\be \label{2.5}
\begin{split}
 & T<\alpha_0<\alpha_1<\dots<\alpha_k , \\
 & T<\beta_1<\beta_2<\dots<\beta_k, \\
 & \alpha_0\in (T,T+H),\ \alpha_r,\beta_r\in (\overset{r}{T},\overset{r}{\wideparen{T+H}}), \\
 & r=1,\dots, k,
\end{split}
\ee
\be \label{2.6}
\begin{split}
 & \alpha_{r+1}-\alpha_r\sim (1-c)\pi(T),\ r=0,1,\dots,k-1, \\
 & \beta_{r+1}-\beta_r\sim (1-c)\pi(T),\ r=1,\dots,k-1,
\end{split}
\ee
where
\bdis
\pi(T)\sim \frac{T}{\ln T},\ T\to\infty
\edis
is the prime-counting function and $c$ is the Euler's constant.
\end{mydef1}

\begin{remark}
Let us notice that the asymptotic behavior of the sets
\be \label{2.7}
\{ \alpha_r\}_{r=0}^k,\ \{\beta_r\}_{r=1}^k
\ee
is as follows: at $T\to\infty$ the points of every set in (\ref{2.7}) recede unboundedly each from other and
all together recede to infinity. Hence, at $T\to\infty$ each set in (\ref{2.7}) looks like one-dimensional
Friedmann-Hubble universe.
\end{remark}

\subsection{}

Let us denote the mean-value of the function
\bdis
\arg\zf,\ t\in [0,T]
\edis
by the symbol
\bdis
\left. \langle \arg\zf\rangle\right|_{[0,T]}.
\edis
Let us mention that the function under consideration has an infinite set of first-order discontinuities. Since
\bdis
\int_0^{\alpha_0(T)}\arg\zf{\rm d}t=\alpha_0(T)\left. \langle \arg\zf\rangle\right|_{[0,\alpha_0(T)]},
\edis
then we obtain from (\ref{2.4}) the following

\begin{mydef41}
\be \label{2.8}
\begin{split}
&\left| \left. \langle \arg\zf\rangle\right|_{[0,\alpha_0(T)]}\right| \sim \\
& \sim \frac{\pi (c_l)^{\frac{1}{2l}}}{\alpha_0(T)}\prod_{r=1}^k
 \left|
\frac
{\zeta\left(\frac 12+i\alpha_r(T,l)\right)}
{\zeta\left(\frac 12+i\beta_r(T,l)\right)}
\right|^{-\frac 1l} , \\
& \alpha_0(T)\in (T,T+H),\ T\to\infty .
\end{split}
\ee
\end{mydef41}

Let us remind that the following Littlewood's estimate (comp. \cite{5}, p. 189)
\bdis
S_1(t)=\mcal{O}(\ln t),\ t\to\infty
\edis
holds true. Hence, we have (see (\ref{1.1}), (\ref{1.2})) the estimate
\be \label{2.9}
\left. \langle \arg\zf\rangle\right|_{[0,T]}=\mcal{O}\left(\frac{\ln T}{T}\right),\ T\to\infty.
\ee

\begin{remark}
Consequently, we have obtained in the direction of the estimate (\ref{2.9}) the explicit asymptotic formula (\ref{2.8}) for the mean-value
\bdis
\left. \langle \arg\zf\rangle\right|_{[0,T]}
\edis
on the infinite subset
\bdis
\{\alpha_0(T)\}, \alpha_0(T)\in (T,T+T^{1/2+\epsilon}),\ T\to\infty.
\edis
\end{remark}

\section{Reduction of the integral in (\ref{2.4})}

\subsection{}

Now, we use the Selberg's $\Omega$-theorem (\ref{1.5}) to transform our formula (\ref{2.4}). It follows from (\ref{1.5}) that
there are two sequences

\bdis
\{ a_n\}_{n=1}^\infty,\ \{ b_n\}_{n=1}^\infty,\ a_n,b_n\to\infty
\edis
such that
\be \label{3.1}
\begin{split}
 & S_1(a_n)>A(\ln a_n)^{1/3}(\ln\ln a_n)^{-10/3}, \\
 & S_1(b_n)<-B(\ln a_n)^{1/3}(\ln\ln a_n)^{-10/3}; \\
 & A,B>0.
\end{split}
\ee
Since
\bdis
S_1(t),\ t>0
\edis
is the continuous function then by (\ref{3.1}) there is (Bolzano-Cauchy) the sequence
\be \label{3.2}
\{\mu_n\}_{n=1}^\infty:\ S_1(\mu_n)=0,\ \mu_n\to\infty,
\ee
where $\mu_n$ is the odd-order root of the equation
\be \label{3.3}
S_1(t)=0,\ t>0 .
\ee

\begin{remark}
We may suppose, of course, that the sequence (\ref{3.2}) is \emph{complete} one in the usual sense, the interval
\bdis
(\mu_n,\mu_{n+1})
\edis
does not contain any other odd-order root of the equation (\ref{3.3}).
\end{remark}

\begin{remark}
There is no need to discuss (for our purpose) the question about even-order roots of the equation (\ref{3.3}).
\end{remark}

Hence, we have: if
\be \label{3.4}
\bar{k}=\bar{k}[\alpha_0(T)]:\ \mu_{\bar{k}}<\alpha_0(T)<\mu_{\bar{k}+1}
\ee
and (of course, see (\ref{2.4}), (\ref{3.3}))
\bdis
S_1[\alpha_0(T)]\not=0,
\edis
then
\be \label{3.5}
\begin{split}
 & S_1[\alpha_0(T)]=\int_0^{\alpha_0(T)}S(t){\rm d}t=\int_0^{\mu_{\bar{k}}}+\int_{\mu_{\bar{k}}}^{\alpha_0(T)}=
 \int_{\mu_{\bar{k}}}^{\alpha_0(T)}S(t){\rm d}t.
\end{split}
\ee
Consequently, we have from (\ref{2.4}) by (\ref{3.4}), (\ref{3.5}) the following

\begin{mydef42}
\be \label{3.6}
\begin{split}
 & \left|\int_{\mu_{\bar{k}}}^{\alpha_0(T)}\arg\zf{\rm d}t\right|\sim \\
 & \sim \pi (c_l)^{\frac{1}{2l}}\prod_{r=1}^k
 \left|
\frac
{\zeta\left(\frac 12+i\alpha_r(T,l)\right)}
{\zeta\left(\frac 12+i\beta_r(T)\right)}
\right|^{-\frac 1l},\ T\to\infty .
\end{split}
\ee
\end{mydef42}

\subsection{}

Next, we obtain from (\ref{3.6}), (comp. (\ref{2.8})), the following
\begin{mydef43}
\be \label{3.7}
\begin{split}
 & \left|\left. \langle \arg\zf\rangle\right|_{[\mu_{\bar{k}},\alpha_0(T)]}\right|\sim \\
 & \sim
 \frac{\pi (c_l)^{\frac{1}{2l}}}{\alpha_0(T)-\mu_{\bar{k}}}\prod_{r=1}^k
 \left|
\frac
{\zeta\left(\frac 12+i\alpha_r(T,l)\right)}
{\zeta\left(\frac 12+i\beta_r(T)\right)}
\right|^{-\frac 1l},
\end{split}
\ee
and, of course, (see (\ref{2.8}), (\ref{3.7}))
\bdis
\begin{split}
 & \left|\left. \langle \arg\zf\rangle\right|_{[\mu_{\bar{k}},\alpha_0(T)]}\right|\sim \\
 & \sim \frac{\alpha_0(T)}{\alpha_0(T)-\mu_{\bar{k}}}
 \left|\left. \langle \arg\zf\rangle\right|_{[0,\alpha_0(T)]}\right|,\ T\to\infty.
\end{split}
\edis
\end{mydef43}

\section{On infinite set of metamorphoses of the Q-system that is generated by the factorization
formula (\ref{2.4})}

\subsection{}

Let us remind the Riemann-Siegel formula
\be \label{4.1}
Z(t)=2\sum_{n\leq\tau(t)}\frac{1}{\sqrt{n}}\cos\{ \vth(t)-t\ln n\}+\mcal{O}(t^{-1/4}),
\ee
where
\bdis
\begin{split}
 & Z(t)=e^{i\vth(t)}\zf,\ \tau(t)=\sqrt{\frac{t}{2\pi}}, \\
 & \vth(t)=-\frac t2\ln\pi+\mbox{Im}\ln\Gamma\left(\frac 14+i\frac t2\right),
\end{split}
\edis
(see \cite{5}, pp. 79, 239). Next, we have introduced (see \cite{2}, (2.1)) the following oscillatory
Q-system (based exactly on the Riemann-Siegel formula (\ref{4.1}))
\be \label{4.2}
\begin{split}
 & G(x_1,\dots,x_k;y_1,\dots,y_k)=\prod_{r=1}^k \left|\frac{Z(x_r)}{Z(y_r)}\right|= \\
 & = \prod_{r=1}^k
 \left|
 \frac
 {\sum_{n\leq\tau(x_r)}\frac{2}{\sqrt{n}}\cos\{\vth(x_r)-x_r\ln n\}+R(x_r)}
 {\sum_{n\leq\tau(y_r)}\frac{2}{\sqrt{n}}\cos\{\vth(y_r)-y_r\ln n\}+R(y_r)}
 \right| , \\
 & (x_1,\dots,x_k)\in M_k^1,\ (y_1,\dots,y_k)\in M_k^2, \\
 & R(t)=\mcal{O}(t^{-1/4}),\ k\leq k_0\in\mbb{N},
\end{split}
\ee
where
\be \label{4.3}
\begin{split}
 & M_k^1=\{ (x_1,\dots,x_k)\in (T_0,+\infty)^k,\ T_0<x_1<\dots<x_k\}, \\
 & M_k^2=\{ (y_1,\dots,y_k)\in (T_0,+\infty)^k,\ T_0<y_1<\dots<y_k\}, \\
 & x_r,y_r\not=\gamma:\ \zeta\left(\frac 12+i\gamma\right)=0,\ r=1,\dots,k.
\end{split}
\ee

\subsection{}

Next, we have obtained (see \cite{2}, (3.1)) the following spectral formula
\be \label{4.4}
\begin{split}
 & Z(t)=2\sum_{n\leq \tau(x_r)}\frac{1}{\sqrt{n}}
 \cos\left\{ t\ln\frac{\tau(x_r)}{n}-\frac{x_r}{2}-\frac{\pi}{8}\right\}+ \\
 & + \mcal{O}(x_r^{-1/4}),\ \tau(x_r)=\sqrt{\frac{x_r}{2\pi}}, \\
 & t\in [x_r,x_r+V],\ V\in (0,\sqrt[4]{x_r}],
\end{split}
\ee
(and similarly for $x_r\longrightarrow y_r$), where
\bdis
T_0<x_r,y_r,\ r=1,\dots,k.
\edis

\begin{remark}
The spectral formula (\ref{4.4}) is, of course, a variant of the Riemann-Siegel formula (\ref{4.1}).
\end{remark}

\begin{remark}
We call the expressions
\be \label{4.5}
\frac{2}{\sqrt{n}}\cos\left\{ t\omega_n(x_r)-\frac{x_r}{2}-\frac{\pi}{8}\right\}, \dots
\ee
as the local Riemann's oscillators with:
\begin{itemize}
 \item[(a)] the amplitudes
 \bdis
 \frac{2}{\sqrt{n}},
 \edis
 \item[(b)] the incohorent local phase constants
 \bdis
 \left\{-\frac{x_r}{2}-\frac{\pi}{8}\right\},\ \left\{-\frac{y_r}{2}-\frac{\pi}{8}\right\},
 \edis
 \item[(c)] the non-synchronized local times
 \bdis
 t=t(x_r)\in [x_r,x_r+V], \dots
 \edis
 \item[(d)] the local spectrum of the cyclic frequencies
 \bdis
 \begin{split}
  & \{\omega_n(x_r)\}_{n\leq\tau(x_r)},\ \omega_n(x_r)=\ln\frac{\tau(x_r)}{n}, \\
  & \{\omega_n(y_r)\}_{n\leq\tau(y_r)},\ \omega_n(y_r)=\ln\frac{\tau(y_r)}{n} .
 \end{split}
 \edis
\end{itemize}
\end{remark}

\begin{remark}
The Q-system (\ref{4.2}) represents a complicated oscillating process generated by oscillations of
big number of the local Riemann's oscillators (\ref{4.5}).
\end{remark}

\subsection{}

Now, in connection with the oscillating Q-system (\ref{4.2}), the following corollary follows from our Theorem

\begin{mydef44}
\be \label{4.6}
\begin{split}
 & \prod_{r=1}^k
 \left|
 \frac
 {\sum_{n\leq\tau(\alpha_r)}\frac{2}{\sqrt{n}}\cos\{\vth(\alpha_r)-\alpha_r\ln n\}+R(\alpha_r)}
 {\sum_{n\leq\tau(\beta_r)}\frac{2}{\sqrt{n}}\cos\{\vth(\beta_r)-\beta_r\ln n\}+R(\beta_r)}
 \right| \sim \\
 & \sim
 \pi^l\sqrt{c_l}
 \left|
 \int_0^{\alpha_0(T)}\arg \zf {\rm d}t
 \right|^{-l},\ T\to\infty.
\end{split}
\ee
\end{mydef44}

\begin{remark}
Hence, we have two resp. one parametric sets of control functions (=Golem's shem) for admissible and fixed
$\epsilon,k$, (see (\ref{2.3})),
\be \label{4.7}
\begin{split}
 & \{\alpha_0(T,l),\alpha_1(T,l),\dots,\alpha_k(T,l)\}, \\
 & \{\beta_1(T),\dots,\beta_k(T)\} , \\
 & T\in (T_0,+\infty),\ l\in\mbb{N},
\end{split}
\ee
of the metamorphoses (\ref{4.6}), (comp. \cite{1},\cite{2}).
\end{remark}

\begin{remark}
The mechanism of the metamorphosis is as follows. Let (comp. (\ref{4.3}), (\ref{4.7}))
\be \label{4.8}
\begin{split}
 & M_k^3=\{\alpha_1(T,l),\dots,\alpha_k(T,l)\}, \\
 & M_k^4=\{\beta_1(T),\dots,\beta_k(T)\},
\end{split}
\ee
where, of course,
\be \label{4.9}
\begin{split}
 & M_k^3\subset M_k^1\subset (T_0,+\infty)^k, \\
 & M_k^4\subset M_k^2\subset (T_0,+\infty)^k.
\end{split}
\ee
Now, if we obtain after random sampling of the points
\bdis
(x_1,\dots,x_k),\ (y_1,\dots,y_k)
\edis
(see the conditions (\ref{4.3})) such that
\be \label{4.10}
\begin{split}
 & (x_1,\dots,x_k)=(\alpha_1(T,l),\dots,\alpha_k(T,l))\in M_k^3, \\
 & (y_1,\dots,y_k)=(\beta_1(T),\dots,\beta_k(T))\in M_k^4,
\end{split}
\ee
(see (\ref{4.8}), (\ref{4.9})), then - at the points (\ref{4.10}) - the Q-system (\ref{4.2}) changes
its old form (=chrysalis) to the new one (=butterfly), and the last ist controlled by the
function $\alpha_0(T)$.
\end{remark}

\subsection{}

Now, we rewrite the formula (\ref{4.6}), (comp. (\ref{3.6})), as follows:
\be \label{4.11}
\begin{split}
 & \left|
 \int_{\mu_{\bar{k}}}^{\alpha_0(T)}\arg\zf{\rm d}t
 \right|\sim \\
 & \sim \pi (c_l)^{\frac{1}{2l}}\prod_{r=1}^k
 \left|
 \frac
 {\sum_{n\leq\tau(\alpha_r)}\frac{2}{\sqrt{n}}\cos\{\vth(\alpha_r)-\alpha_r\ln n\}+R(\alpha_r)}
 {\sum_{n\leq\tau(\beta_r)}\frac{2}{\sqrt{n}}\cos\{\vth(\beta_r)-\beta_r\ln n\}+R(\beta_r)}
 \right|^{-\frac 1l}.
\end{split}
\ee

\begin{remark}
The formula (\ref{4.11}) expresses the metamorphosis in the reverse direction. We describe the mechanism of
this as follows: we begin with the integral
\bdis
\left|
\int_0^w \arg\zf{\rm d}t
\right|
\edis
that is the Aaron staff,
\bdis
\longrightarrow
\left|
\int_{\mu_{\bar{k}}}^{\alpha_0(T)} \arg\zf{\rm d}t
\right|
\edis
that is the bud of the Aaron staff corresponding to $w=\alpha_0(T)$,
\bdis
\sim
\pi (c_l)^{\frac{1}{2l}}\prod_{r=1}^k
 \left|
 \frac
 {\sum_{n\leq\tau(\alpha_r)}\frac{2}{\sqrt{n}}\cos\{\vth(\alpha_r)-\alpha_r\ln n\}+R(\alpha_r)}
 {\sum_{n\leq\tau(\beta_r)}\frac{2}{\sqrt{n}}\cos\{\vth(\beta_r)-\beta_r\ln n\}+R(\beta_r)}
 \right|^{-\frac 1l}
\edis
already metamorphosed one into almonds ripened, (motivation: Chumash, Bamidbar, 17:23).
\end{remark}

\thanks{I would like to thank Michal Demetrian for his help with electronic version of this paper.}

\end{document}